\documentclass{article}
\usepackage{verbatim}
\usepackage{amsmath}
\usepackage[mathscr]{eucal}
\usepackage{amsthm}
\usepackage{amssymb}
\usepackage[english,german]{babel}
\theoremstyle{plain}
\newtheorem{Theo}{Theorem}[section]
\newtheorem{Prop}[Theo]{Proposition}
\newtheorem{Lem}[Theo]{Lemma}
\newtheorem{Cor}[Theo]{Corollary}
\newtheorem{Main}{Theorem}
\theoremstyle{definition}
\newtheorem{Def}[Theo]{Definition}
\theoremstyle{remark}

\begin{document}

\selectlanguage{english}
\title{The $q$-analogue of the alternating groups and its representations}
\author{Hideo Mitsuhashi \\
Fujisawa Advanced Vocational Skill Training School \\ 
290--2 Kawana, Fujisawa--shi, Kanagawa--ken 251--0015, Japan}

%\address{Fujisawa Vocational Training School, 
%290-2 Kawana, Fjisawa-shi, Kanagawa-ken 251-0015, Japan} 

\date{}
\maketitle

\section{Introduction}\par
Frobenius began the study of representation theory and character 
theory of the symmetric groups $\frak{S}_{n}$ at the turn of the century [4]. 
In this study, He showed that there is one irreducible representation 
${\pi}_{\lambda}$ of 
$\frak{S}_{n}$ corresponding to each partition $\lambda$ of $n$ and this 
correspondence gives all of irreducible representations of $\frak{S}_{n}$ up to 
equivalency. After this work, he also began the study of representation 
theory and character theory of the alternating groups $\frak{A}_{n}$ 
at the opening of the 20th century [5]. 
In this study, he showed that there are branching rules 
when we restrict irreducible 
representation of the $\frak{S}_{n}$ to $\frak{A}_{n}$. 
We denote by $\Tilde{\pi}_{\lambda}$ the restriction of $\pi_\lambda$ to $\frak{A}_{n}$ . 
If $\lambda$ is  
self-conjugate(that is symmetric figure as Young diagram along its 
diagonal), then $\Tilde{\pi}_{\lambda}$ is decomposed into the direct sum of 
two irreducible representations whose degrees are coincide. But these two 
representations are mutually inequivalent. On the other hand, if $\lambda$ is non 
self-conjugate, then $\Tilde{\pi}_{\lambda}$ is also the irreducible and 
$\Tilde{\pi}_{\lambda}$ is equivalent to $\Tilde{\pi}_{^t\lambda}$ where 
${^t\lambda}$ is transpose of ${\lambda}$. These representations consist all of 
irreducible representations of $\frak{A}_{n}$ up to equivalency(Theorem 3.4). 
The concrete construction of irreducible representations of $\frak{S}_{n}$ were 
given by Young. He gave irreducible representations of $\frak{S}_{n}$ 
in two ways: ${\lq}$seminormal form', and ${\lq}$orthogonal form' (For detail, see [2]). \par
The Hecke algebra $\mathcal{H}_{n}(q)$ of type $A$ is considered as the $q$-analogue 
of the group ring $\mathbb{C}[\frak{S}_{n}]$ of the symmetric group 
$\frak{S}_{n}$, and its irreducible 
representations are constructed as well as in the case of $\frak{S}_{n}$ under some 
assumptions about the value of the parameter $q$. 
The seminormal form for $\mathcal{H}_{n}(q)$ was given by Hoefsmit [9], and the 
orthogonal form was given by Wenzl [15]. 
It was not seem to be seen about the $q$-analogue of the group ring 
$\mathbb{C}[\frak{A}_{n}]$ of the alternating group $\frak{A}_{n}$ so far. 
In this context, we defined a subalgebra $\frak{A}_{n}(q)$ of $\mathcal{H}_{n}(q)$ 
whose dimension as vector space over $\mathbb{C}$ 
is just a half of that of $\mathcal{H}_{n}(q)$. 
The group ring $\mathbb{C}[\frak{A}_{n}]$ which is generated by following elements, 
\[ \{s_1s_{i+1} \mid 1 \leq i \leq n-2 \} \]
where $s_i$'s are simple permutations, has the dimension just a half of $\frak{S}_{n}$ 
over $\mathbb{C}$. 
On the other hand, the subalgebra of $\mathcal{H}_{n}(q)$ generated by following elements, 
\[ \{g_1g_{i+1} \mid 1 \leq i \leq n-2 \} \]
where $g_i$'s are generators of $\mathcal{H}_{n}(q)$, has the dimension more than 
a half of that of $\mathcal{H}_{n}(q)$ because of quadratic relations, 
\[ g_i^2=(q-1)g_i+q \qquad \text{for $i=1,2,\hdots ,n-1$} \]
of $\mathcal{H}_{n}(q)$. So instead of $g_i$, we take another generators $f_i$ 
of $\mathcal{H}_{n}(q)$ as follows.
\[ f_i=\dfrac{2g_i-(q-1)}{q+1} \qquad \text{for $i=1,2,\hdots ,n-1$} \]
Using these generators, we define the subalgebra $\frak{A}_{n}(q)$ of $\mathcal{H}_{n}(q)$ 
generated by following elements.
\[ \{f_1f_{i+1} \mid 1 \leq i \leq n-2 \} \]
Just then, $\frak{A}_{n}(q)$ has the dimension 
just a half of that of $\mathcal{H}_{n}(q)$(Theorem 5.3).
Moreover, when we take the limit $q \rightarrow 1$, 
$\frak{A}_{n}(q)$ goes to $\mathbb{C}[\frak{A}_{n}]$ and the generator 
$f_1f_{i+1}$ goes to $s_1s_{i+1}$. Therefore, we consider $\frak{A}_{n}(q)$ as the $q$-analogue 
of the alternating group $\frak{A}_{n}$. Furthermore, we show a presentation of 
$\frak{A}_{n}(q)$(Theorem 5.5). In this presentation, 
quadratic relations are same as in the case of 
$\frak{A}_{n}$, but cubic relations are slightly different. At the limit 
$q \rightarrow 1$, the defining relations of $\frak{A}_{n}(q)$ coincide with 
those of $\frak{A}_{n}$. Thus we can define the $q$-analogue of the alternating groups 
$\frak{A}_{n}$ as follows.
\begin{Def}
Assume that $q\neq -1$ and $n>2$. 
$\frak{A}_{n}(q)$ is the algebra over $\mathbb{C}$ defined by the generators 
$y_1,y_2,\hdots ,y_{n-2}$ and the following relations.
\begin{enumerate}
\rm
\item
\it
$y_1^3=-\Big{(}\dfrac{q-1}{q+1}\Big{)}^2(y_1^2-y_1)+1$
\rm
\item
\it
$y_i^2=1 \qquad \text{for $i>1$}$
\rm
\item
\it
$(y_{i-1}y_i)^3=-\Big{(}\dfrac{q-1}{q+1}\Big{)}^2\Big{\{}(y_{i-1}y_i)^2-y_{i-1}y_i\Big{\}}+1 
\qquad \text{for $i=2,3,\hdots ,n-2$}$
\rm
\item
\it
$(y_iy_j)^2=1 \qquad \text{whenever $|i-j|>1$}$
\end{enumerate}
For $n=2$, we define $\frak{A}_{2}(q)$ the algebra over $\mathbb{C}$ 
generated by only the unit element.
\end{Def}
We also analyze the representations of $\frak{A}_{n}(q)$. Using the same notation as in 
the case of $\frak{S}_{n}$(no confusion appears in this context), 
for the irreducible representation 
${\pi}_{\lambda}$ of $\mathcal{H}_{n}(q)$ corresponding to the partition $\lambda$ of $n$ 
(we identify $\lambda$ with a certain Young tableau), we denote by $\Tilde{\pi}_{\lambda}$ 
the restriction of $\pi_\lambda$ to $\frak{A}_{n}(q)$. 
Except of some $\lq$bad' values of $q$, branching rules for $\frak{S}_{n}\rightarrow 
\frak{A}_{n}$ are valid to $\mathcal{H}_{n}(q)\rightarrow \frak{A}_{n}(q)$. 
More precisely, if $\lambda$ is non self-conjugate then $\Tilde{\pi}_{\lambda}$ is irreducible, 
on the other hand, if $\lambda$ is self-conjugate then $\Tilde{\pi}_{\lambda}$ is decomposed 
into two irreducible 
components $\Tilde{\pi}_{\lambda}^+$ and $\Tilde{\pi}_{\lambda}^-$. The following statement is 
the main result for representations of $\frak{A}_{n}(q)$.
\renewcommand{\theMain}{6.5}
\begin{Main}
Let $q$ be a complex number such that $q\neq 0$ and $q$ is not a 
$k$-th root of unity with $1\leq k \leq n$. 
Let $\lambda_1,^t\lambda_1,\lambda_2,^t\lambda_2,\hdots,
\lambda_p,^t\lambda_p$ be non self-conjugate 
Young diagrams and $\lambda_{p+1},\lambda_{p+2},\hdots,\lambda_{p+q}$ 
be self-conjugate Young diagrams within $\Lambda_n$. 
Then representations 
$\Tilde{\pi}_{\lambda_1}$,$\Tilde{\pi}_{\lambda_2}$,
$\hdots$,$\Tilde{\pi}_{\lambda_p}$,$\Tilde{\pi}_{\lambda_{p+1}}^+$,
$\Tilde{\pi}_{\lambda_{p+1}}^-$,$\hdots,\Tilde{\pi}_{\lambda_{p+q}}^+$,
$\Tilde{\pi}_{\lambda_{p+q}}^-$ are irreducible and not equivalent 
each other. 
These representations consist of all equivalent classes of 
irreducible representations of $\frak{A}_{n}(q)$. Hence, 
$\frak{A}_{n}(q)$ is semisimple.
\end{Main}
We show this result by using relations of matrix elements of orthogonal representations 
of $\mathcal{H}_{n}(q)$. Thus, we get all of irreducible representations of 
$\frak{A}_{n}(q)$ up to equivalence and show the semisimplicity of $\frak{A}_{n}(q)$ 
except of some $\lq$bad' values of $q$.

\section{Orthogonal representations of the symmetric groups}\par
In this section, we give a review of the concrete construction of 
an irreducible representation corresponding to each Young diagram, 
It is called ${\lq}$Orthogonal representation', which was given by A.Young 
in 1932. \par
Let $\frak{S}_{n}$ be the symmetric group, consisting of all 
permutations of $n$ things. Especially, $\frak{S}_{n}$ 
permutes $n$ numbers $ \{ 1,2,\hdots,n \} $. 
$\frak{S}_{n}$ is generated by elements $s_i$($i=1,2,\hdots,n-1$) 
where $s_i=(i,i+1)$ is a simple permutation, and have defining relations as follows. 
\begin{enumerate}
\item
$s_i^2=1 \qquad \text{for $i=1,2,\hdots ,n-1$}$
\item
$s_is_{i+1}s_i=s_{i+1}s_is_{i+1} \qquad \text{for $i=1,2,\hdots ,n-2$}$
\item
$s_is_j=s_js_i \qquad \text{whenever $|i-j|\geq 2$}$
\end{enumerate}
We notice that $\left| \frak{S}_{n} 
\right| = n!$. 
For $\sigma \in \frak{S}_{n}$
\[ \sigma = s_{i_1}s_{i_2} \hdots s_{i_r} \]
is an expression of $\sigma$. When $\sigma$ cannot be expressed by 
any form
\[ \sigma = s_{j_1}s_{j_2} \hdots s_{j_t} \]
with $t<r$, we call $l(\sigma) = r$ the length of $\sigma$ 
and the expression
\[ \sigma = s_{i_1}s_{i_2} \hdots s_{i_r} \]
is called the reduced expression of $\sigma$.\par
Consider a permutation
\[ \sigma = \begin{pmatrix}
1 & 2 & \hdots & n \\
\sigma_1 & \sigma_2 & \hdots & \sigma_n
\end{pmatrix}
\]
and set the number of inversions of $\sigma$ as follows.
\[ d_i = \sharp \{ j \mid j>k \qquad \text{where $\sigma_k=i$ 
and $\sigma_j<i$} \} \]
Then, $\sigma$ has a reduced expression
\begin{equation*}
\begin{split}
\sigma &= \hdots (s_{i-1}s_{i-2}\hdots s_{i-{d_i}}) \\
& \quad \hdots (s_{n-2}s_{n-3}\hdots s_{n-1-{d_{n-1}}}) \\
& \quad (s_{n-1}s_{n-2}\hdots s_{n-{d_n}})
\end{split}
\end{equation*}
where the $i$-th contribution is only to be included if $d_i \geq 1$. 
The reduced expression described above is called the normal reduced 
expression. For the sake of convenience, we denote $s_{i-1}s_{i-2}\hdots s_{i-{d_i}}$ 
by $U_{i,d_i}$.
\begin{Prop}
The map
\[ f:[0,1]\times[0,2]\times\hdots\times[0,n-1]\longrightarrow
\frak{S}_{n} \]
defined by
\[ f(d_2,d_3,\hdots,d_n)=
U_{2,d_2}U_{3,d_3} \hdots U_{i,d_i} \hdots U_{n,d_n} \]
(where the $i$-th contribution is only to be included if $d_i \geq 1$, 
and $f(0,0,\hdots,0)$ is the unit element in $\frak{S}_{n}$)\\
is a bijection.
\end{Prop}
\begin{Cor}
Let $\mathbb{C}[\frak{S}_{n}]$ be the group ring of 
$\frak{S}_{n}$. Then $\mathbb{C}[\frak{S}_{n}]$ has following basis 
as a vector space over $\mathbb{C}$.
\[ \{U_{2,d_2}U_{3,d_3} \hdots U_{n,d_n}\mid 
(d_2,d_3,\hdots,d_n)\in [0,1]\times[0,2]\times\hdots\times[0,n-1]
\} \]
\end{Cor}
Let $\Lambda_n$ be the 
set of all Young diagrams with $n$ boxes. We mean by a Young diagram
$\lambda = [\lambda_1,\lambda_2,\ldots ] \in \Lambda_n$ an array of
$n$ boxes with $\lambda_1$ boxes in the first row, $\lambda_2$ boxes 
in the second row, and so on with $\lambda_1 \geq \lambda_2 \geq\hdots$. 
For $\lambda \in \Lambda_n$, the number of rows in which there are nonzero 
boxes is called the depth of $\lambda$ and we denote by $l(\lambda)$. 
We immediately observe that $\lambda_1+\lambda_2+\hdots\lambda_{l(\lambda)}=n$.
Formally, the Young diagram $\lambda \in \Lambda_n$ is defined 
by the subset of $\mathbb{N} \times \mathbb{N}$ as follows.
\[ \{ (i,j)\in \mathbb{N} \times \mathbb{N} \mid 
1 \leq j \leq \lambda_i, 1 \leq i \leq l(\lambda) \} \]
The element $(i,j)$ is identified with the box at $i$-th row 
and $j$-th column of $\lambda$. We shall write $\mu < \lambda$ if $\mu$ is 
obtained from $\lambda$ by removing appropriate boxes.\par
As well known fact, for each Young diagram $\lambda \in \Lambda_n$, 
there is an irreducible representation of $\frak{S}_{n}$ over
the complex field $\mathbb{C}$ up to equivalence, and all irreducible 
representations of $\frak{S}_{n}$ over $\mathbb{C}$ are parameterized by 
elements of $\Lambda_n$.\par
Let $\operatorname{Tab}(\lambda)$ be the set of all tableaux belonging to $\lambda$ 
and let $\operatorname{STab}(\lambda)$ be the set of all standard tableaux 
belonging to $\lambda$. For $T \in \operatorname{Tab}(\lambda)$, We mean by $T(i,j)$ 
the number on the box at $i$-th row and $j$-th column of $T$.\par
$\frak{S}_{n}$ acts on $\operatorname{Tab}(\lambda)$ as 
permutations of numbers on boxes. For $\sigma \in \frak{S}_{n}$ 
and $T \in \operatorname{Tab}(\lambda)$, $\sigma$ acts on $T$ as follows.
\[ (\sigma \cdot T)(i,j) = \sigma (T(i,j)) \]
We claim that this action is simply transitive.
\begin{Def}
Let $T$ be a tableau of shape $\lambda$, and $(i,j)$ the box at $i$-th row and 
$j$-th column of $T$. Then the class of $(i,j)$ is the value $j-i$. 
For $k = T(i,j)$, we mean by $\alpha_{k,T}$ the class of $(i,j)$.
\end{Def}
From Definition 2.3, we have $\alpha_{1,T} = 0$ for each standard tableau $T$.
\begin{Def}
Let $k$,$l$ be numbers on boxes in a tableau $T$, 
then the axial distance $d_{T,k,l}$ from $k$ to $l$ in $T$
is defined as follows.
\[ d_{T,k,l} = \alpha_{k,T}-\alpha_{l,T} \]
\end{Def}
It follows immediately from Definition 2.4 that
\[ d_{T,k,l} = -d_{T,l,k} \]
\begin{Lem}
If $T \in \operatorname{STab}(\lambda)\text{($\lambda \in \Lambda_n$)}$, 
then $d_{T,k,k+1}$ is nonzero for each $k$ with $1 \leq k \leq n-1$.
\end{Lem}
\begin{proof}
If $(i,j)$ and $(i',j')$ have a same class,
then there exists an integer $m$ such that $i=i'+m,j=j'+m$. Hence, boxes which have a 
same class are located on the line parallel to the diagonal of 
the Young diagram.  For standard tableaux, numbers labeled on 
boxes strictly increase across each row and each column. 
Therefore, the box on which $k$ is labeled and the box on which $k+1$ 
is labeled cannot be located on the same line parallel to diagonal.
\end{proof}
Let $V_\lambda$ be the complex vector space with basis 
$\{v_T \mid T \in \operatorname{STab}(\lambda) \}$.
We define the linear map $\pi_\lambda(s_i)$ for elements 
$s_i = (i,i+1) \in \frak{S}_{n}$ on $V_\lambda$ in three ways 
according to the relation of location between the box on which 
$i$ is labeled and the box on which $i+1$ is labeled in $T$.
\begin{enumerate}
\item
If $i$ and $i+1$ appear in the same row, 
$\pi_\lambda(s_i)v_T=v_T$.
\item
If $i$ and $i+1$ appear in the same column, 
$\pi_\lambda(s_i)v_T=-v_T$.
\item
Otherwise, we observe that $s_i \cdot T$ is also in 
$\operatorname{STab}(\lambda)$, and $\pi_\lambda(s_i)$ acts on the 
subspace $\mathbb{C}v_T \oplus \mathbb{C}v_{s_i\cdot T}$ of $V_\lambda$ 
defined by the matrix,	
\[
\begin{bmatrix}
-\eta & \sqrt{1-\eta^2} \\
\sqrt{1-\eta^2} & \eta 
\end{bmatrix}
\]
where $\eta = (d_{T,i,i+1})^{-1}$.
\end{enumerate}
\begin{Theo}
For $\lambda \in \Lambda_n$ 
and $\mu \in \Lambda_n$ such that $\lambda \neq \mu$, we have  
$\pi_\lambda \ncong \pi_\mu$. Moreover, by the definition above, 
$\{\pi_\lambda \mid 
\lambda \in \Lambda_n \}$ is the set of all equivalent classes of 
irreducible representations of $\frak{S}_{n}$.
\end{Theo}
\begin{proof}
See for example, [6] or [15].
\end{proof}

\section{Representations of the alternating groups}\par
In this section, we review some facts about irreducible 
representations of the alternating groups. Every irreducible representation 
of the alternating groups is realized with the restriction of that of the 
symmetric groups. Furthermore, there is a beautiful branching rule for the 
restriction of irreducible representations of the symmetric groups to the 
alternating groups, which was discovered by Frobenius [5]. \par
Let $\frak{A}_{n}$ be the alternating group, consisting of all even 
permutations. 
\[ \frak{A}_{n}=\{ \sigma \in \frak{S}_{n} \mid 
\text{$l(\sigma)$ is even} \} \]
As the fact, we can say that $\frak{A}_{n}$ is generated by 
following elements.
\[ \frak{A}_{n}=\langle s_1s_{i+1} \mid 1 \leq i \leq n-2 \rangle \]
We notice that $\left| \frak{A}_{n} \right| = n!/2 
= \left| \frak{S}_{n} \right| /2$.\par
From Corollary 2.2 and the definition of $\frak{A}_{n}$, 
the following statement holds.
\begin{Prop}
Let $\mathbb{C}[\frak{A}_{n}]$ be the group ring of 
$\frak{A}_{n}$. Then $\mathbb{C}[\frak{A}_{n}]$ has 
following basis as a vector space over $\mathbb{C}$.
\begin{equation*}
\begin{split}
 \{\sigma=U_{2,d_2}U_{3,d_3} \hdots U_{n,d_n}\mid & (d_1,d_2,\hdots,d_{n-1}) \\ 
& \in [0,1]\times[0,2]\times\hdots\times[0,n-1] \quad \text{$l(\sigma)$ is even} \}
\end{split}
\end{equation*}
\end{Prop}
For irreducible representations $\pi_\lambda$ of $\frak{S}_{n}$, 
we denote by $\Tilde{\pi}_\lambda$ the restriction of $\pi_\lambda$ to 
$\frak{A}_{n}$.
For $\lambda \in \Lambda_n$ (resp. $T \in \operatorname{STab}(\lambda)$), we mean by 
$^t\lambda$ (resp. $^tT$) its transpose.
When $^t\lambda = \lambda$, we say $\lambda$ self-conjugate.
\begin{Prop}
For every non self-conjugate $\lambda \in \Lambda_n$, 
$\Tilde{\pi}_\lambda$ is an irreducible representation of 
$\frak{A}_{n}$ and
$\Tilde{\pi}_\lambda \cong \Tilde{\pi}_{^t\lambda}$.
\end{Prop}
\begin{Prop}
For every self-conjugate $\lambda \in \Lambda_n$, 
$\Tilde{\pi}_\lambda$ is the direct sum of two irreducible representations 
$\Tilde{\pi}^+_\lambda$ and $\Tilde{\pi}^-_\lambda$ 
with $\operatorname{deg}(\Tilde{\pi}^+_\lambda) = \operatorname{deg}(\Tilde{\pi}^-_\lambda)$.
\end{Prop}
\begin{Theo}
Let $\lambda_1$,$^t\lambda_1$,$\lambda_2$,$^t\lambda_2$,$\hdots$,
$\lambda_p$,$^t\lambda_p$ be non self-conjugate 
Young diagrams and $\lambda_{p+1}$,$\lambda_{p+2}$,$\hdots$, $\lambda_{p+q}$ 
be self-conjugate Young diagrams within $\Lambda_n$. Then irreducible 
representations $\Tilde{\pi}_{\lambda_1}$,$\Tilde{\pi}_{\lambda_2}$,
$\hdots$,$\Tilde{\pi}_{\lambda_p}$,$\Tilde{\pi}_{\lambda_{p+1}}^+$,
$\Tilde{\pi}_{\lambda_{p+1}}^-$,$\hdots$,$\Tilde{\pi}_{\lambda_{p+q}}^+$,
$\Tilde{\pi}_{\lambda_{p+q}}^-$ are not equivalent each other and 
these consist of all equivalent classes of irreducible 
representations of $\frak{A}_{n}$.
\end{Theo}
\begin{proof}
See, for example [2] or [5].
\end{proof}
It is well known that there is a presentation of $\frak{A}_{n}$ as follows. 
\begin{Prop}
For $n>2$, the generators $x_1,x_2,\hdots ,x_{n-2}$ 
and the following relations define a presentation of $\frak{A}_{n}$.
\begin{enumerate}
\rm
\item
\it
$x_1^3=x_2^2= \hdots =x_{n-2}^2=1$
\rm
\item
\it
$(x_{i-1}x_i)^3=1 \qquad \text{for $i=2,3,\hdots ,n-2$}$
\rm
\item
\it
$(x_ix_j)^2=1 \qquad \text{whenever $|i-j|>1$}$
\end{enumerate}
\end{Prop}

\section{Orthogonal representations of the Hecke algebra 
$\mathcal{H}_{n}(q)$}\par
Like the case of $\frak{S}_{n}$, all equivalent classes of 
irreducible representations of the Hecke algebra $\mathcal{H}_{n}(q)$ 
are parametrized by the elements in $\Lambda_n$. 
In this section, we will review how to 
construct an irreducible representation of $\mathcal{H}_{n}(q)$ 
corresponding to $\lambda \in \Lambda_n$, which was taken from 
Wenzl's paper [15].\par
Since the way of construction is similar to the case of 
$\frak{S}_{n}$, and we do not discuss about representations of 
$\frak{S}_{n}$ any more, 
we will use the same notation as in the case of $\frak{S}_{n}$.\par
Let $q$ be a complex number and we will mean by the Hecke algebra 
$\mathcal{H}_{n}(q)$ of type $A_{n-1}$ the algebra over 
$\mathbb{C}$ with generators $1,g_1,g_2,\hdots ,g_{n-1}$ with 
the defining relations
\begin{enumerate}
\item
$g_i^2=(q-1)g_i+q \qquad \text{for $i=1,2,\hdots ,n-1$}$
\item
$g_ig_{i+1}g_i=g_{i+1}g_ig_{i+1} \qquad \text{for $i=1,2,\hdots ,n-2$}$
\item
$g_ig_j=g_jg_i \qquad \text{whenever $|i-j|\geq 2$}$
\end{enumerate}
When $q$ is not a root of unity, it is known that $\mathcal{H}_{n}(q)
\cong \mathbb{C}[\frak{S}_{n}]$.\par
If $g_{i_1}g_{i_2}\hdots g_{i_k}$ cannot be expressed by 
any linear combination whose terms are products of at most $k-1$ generators, 
then $g_{i_1}g_{i_2}\hdots g_{i_k}$ is called the reduced expression.
Using the same notation as in the case of $\mathbb{C}[\frak{S}_{n}]$, 
we denote $g_{i-1}g_{i-2}\hdots g_{i-{d_i}}$ 
by $U_{i,d_i}$. As well known fact, we can show the following.
\begin{Prop}
$\mathcal{H}_{n}(q)$ has following basis as a vector space over $\mathbb{C}$.
\[ \{U_{2,d_2}U_{3,d_3}\hdots U_{i,d_i} \hdots U_{n,d_n}\mid 
(d_2,d_3,\hdots,d_n)\in [0,1]\times[0,2]\times\hdots\times[0,n-1]
\} \]
(where the $i$-th contribution is only to be included if $d_i \geq 1$, 
and $U_{2,0}U_{3,0}\hdots U_{n,0}$ is the unit element in $\mathcal{H}_{n}(q)$)
\end{Prop}
Let $V_\lambda$ be the complex vector space with basis 
$\{v_T \mid T \in \operatorname{STab}(\lambda) \}$ for $\lambda \in \Lambda_n$.
We define the linear map $\pi_\lambda(g_i)$ for elements 
$g_i \in \mathcal{H}_{n}(q)$ on $V_\lambda$ in three ways 
according to the relation of location between the box on which 
$i$ is labeled and the box on which $i+1$ is labeled in $T$.
\begin{enumerate}
\item
If $i$ and $i+1$ appear in the same row, 
$\pi_\lambda(g_i)v_T=qv_T$.
\item
If $i$ and $i+1$ appear in the same column, 
$\pi_\lambda(g_i)v_T=-v_T$.
\item
Otherwise, we observe that $s_i \cdot T$ is also in 
$\operatorname{STab}(\lambda)$, and $\pi_\lambda(g_i)$ acts on the 
subspace $\mathbb{C}v_T \oplus \mathbb{C}v_{s_i\cdot T}$ of $V_\lambda$ 
defined by the matrix, 
\[
\dfrac{1}{1-q^d}
\begin{bmatrix}
	q^d(1-q) & \sqrt{q(1-q^{d-1})(1-q^{d+1})} \\
	&	\\
	\sqrt{q(1-q^{d-1})(1-q^{d+1})} & -(1-q) 
\end{bmatrix}
\]
where $d = d_{T,i,i+1}$ is the axial distance from $i$ to $i+1$ in $T$.
\end{enumerate}
When $q\neq 0$ and $q$ is not a $k$-th root of unity with $1\leq k \leq n-1$, 
this matrix is always well-defined because of Lemma 2.5.
\begin{Theo}
For $\lambda \in \Lambda_n$ 
and $\mu \in \Lambda_n$ such that $\lambda \neq \mu$, we have 
$\pi_\lambda \ncong \pi_\mu$. Moreover, by the definition above, 
$\{\pi_\lambda \mid 
\lambda \in \Lambda_n \}$ is all equivalent classes of 
irreducible representations of $\mathcal{H}_{n}(q)$.
\end{Theo}
To see the relation of representation matrix between $\pi_\lambda$ 
and $\pi_{^t\lambda}$, we show the following lemma.
\begin{Lem}
For any $T \in \operatorname{STab}(\lambda)$ with $\lambda \in \Lambda_n$, 
following holds.
\[ d_{T,i,j} = -d_{^tT,i,j} \]
\end{Lem}
Therefore, we obtain the following immediately.
\begin{Prop}
The following relations between $\pi_\lambda$ and $\pi_{^t\lambda}$ hold for 
every $\lambda \in \Lambda_n$.
\begin{enumerate}
\rm
\item
\it
If $\pi_\lambda(g_i)v_T=qv_T$, 
then $\pi_{^t\lambda}(g_i)v_{^tT}=-v_{^tT}$
\rm
\item
\it
If $\pi_\lambda(g_i)v_T=-v_T$, 
then $\pi_{^t\lambda}(g_i)v_{^tT}=qv_{^tT}$
\rm
\item
\it
Otherwise, we observe that $s_i \cdot ^tT$ is also in 
$\operatorname{STab}(^t\lambda)$, and $\pi_{^t\lambda}(g_i)$ leaves invariant the 
subspace $\mathbb{C}v_{^tT}\oplus \mathbb{C}v_{s_i \cdot ^tT}$ of $V_{^t\lambda}$.
On the subspace $\mathbb{C}v_{^tT} \oplus \mathbb{C}v_{s_i \cdot ^tT}$, 
$\pi_{^t\lambda}(g_i)$ is described by the matrix
\[
\dfrac{1}{1-q^d}
\begin{bmatrix}
	-(1-q) & -\sqrt{q(1-q^{d-1})(1-q^{d+1})} \\
	&	\\
	-\sqrt{q(1-q^{d-1})(1-q^{d+1})} & q^d(1-q) 
\end{bmatrix}
\]
where $d = d_{T,i,i+1}$ is the axial distance from $i$ to $i+1$ in $T$.
\end{enumerate}
\end{Prop}
We introduce another presentation of $\mathcal{H}_{n}(q)$. 
Assume $q \neq -1$ and set
\[
f_i=\dfrac{2g_i-(q-1)}{q+1} \qquad \text{for $i=1,2,\hdots n-1$}
\] 
\begin{Prop}
These generate $\mathcal{H}_{n}(q)$ and constitute with the relations
\begin{enumerate}
\rm
\item
\it
$f_i^2=1 \qquad \text{for $i=1,2,\hdots ,n-1$}$
\rm
\item
\it
$f_if_{i+1}f_i=f_{i+1}f_if_{i+1}+\Big{(}\dfrac{q-1}{q+1}\Big{)}^2(f_{i+1}-f_i)
 \qquad \text{for $i=1,2,\hdots ,n-1$}$
\rm
\item
\it
$f_if_j=f_jf_i \qquad \text{whenever $|i-j|\geq 2$}$
\end{enumerate}
a presentation of $\mathcal{H}_{n}(q)$.
\end{Prop}
\begin{proof}
A straightforward computation.
\end{proof}
We denote $f_{i-1}f_{i-2} \hdots f_{i-d_i}$ with $i=2,3,\hdots n$ and $d_i \in [0,i-1]$ 
by ${U'}_{i,d_i}$.
In the same way as Proposition 4.1, we can show the following.
\begin{Prop}
$\mathcal{H}_{n}(q)$ has following basis as a vector space over $\mathbb{C}$.
\begin{equation*}
\begin{split}
\{{U'}_{2,d_2}{U'}_{3,d_3} \hdots {U'}_{i,d_i} \hdots {U'}_{n,d_n} \mid 
& (d_2,d_3,\hdots,d_n) \\ & \in [0,1]\times[0,2]\times\hdots\times[0,n-1]
\}
\end{split}
\end{equation*}
(where the $i$-th contribution is only to be included if $d_i \geq 1$, 
and ${U'}_{2,0}{U'}_{3,0}\hdots {U'}_{n,0}$ is the unit element in $\mathcal{H}_{n}(q)$)
\end{Prop}
\begin{proof}
Expanding ${U'}_{2,d_2}$${U'}_{3,d_3}$$\hdots$${U'}_{i,d_i}$$\hdots$${U'}_{n,d_n}$ with $g_i$'s, 
the term which has the longest length is the nonzero scalar 
multiplication of $U_{2,d_2}$$U_{3,d_3}$$\hdots$$U_{i,d_i}$$\hdots$$U_{n,d_n}$ 
from the definition of $f_i$'s. Hence 
${U'}_{2,d_2}$${U'}_{3,d_3}$$\hdots$ \\
${U'}_{i,d_i}$$\hdots$${U'}_{n,d_n}$'s are 
linearly independent and constitute a basis of $\mathcal{H}_{n}(q)$.
\end{proof}
From the definition of $\pi_{\lambda}$ and Proposition 4.4, we observe the following.
\begin{Prop}
The following relations between $\pi_\lambda$ and $\pi_{^t\lambda}$ hold for 
every $\lambda \in \Lambda_n$.
\begin{enumerate}
\rm
\item
\it
$\pi_\lambda(f_i)v_T=v_T$ and $\pi_{^t\lambda}(f_i)v_{^tT}=-v_{^tT}$ 
if $i$ and $i+1$ appear in the same row of $T$.
\rm
\item
\it
$\pi_\lambda(f_i)v_T=-v_T$ and $\pi_{^t\lambda}(f_i)v_{^tT}=v_{^tT}$ 
if $i$ and $i+1$ appear in the same column of $T$.
\rm
\item
\it
Otherwise, the representation matrix of 
$\pi_\lambda(f_i)$ acting on subspace of 
$V_\lambda$ spanned by $(v_T,v_{s_i\cdot T})$, and 
$\pi_{^t\lambda}(f_i)$ acting on subspace of 
$V_{^t\lambda}$ spanned by $(v_{^tT},v_{s_i\cdot{^tT}})$ is given by
\[ 
\pi_\lambda(f_i) \sim
\begin{bmatrix}
	A & B \\
	B & -A \\
\end{bmatrix} \qquad 
\pi_{^t\lambda}(f_i) \sim 
\begin{bmatrix}
	-A & -B \\
	-B & A \\
\end{bmatrix}
\]
where
\begin{align*}
A &= \dfrac{(1-q)(1+q^d)}{(1+q)(1-q^d)} \\
B &= \dfrac{2\sqrt{q(1-q^{d-1})(1-q^{d+1})}}{(1+q)(1-q^d)} \\
d &= d_{T,i,i+1} \\
\end{align*}
\end{enumerate}
\end{Prop}

\section{The subalgebra $\frak{A}_{n}(q)$ of $\mathcal{H}_{n}(q)$}\par
From Proposition 4.5, we observe that 
parity of length of each expression by $f_i$'s is preserved 
even if we choose another expression. To be exact, if an expression of even (resp. odd) 
length is expressed in a linear combination of other expressions, then the length of 
each term is even (resp. odd). So we can consider the subalgebra $\frak{A}_{n}(q)$ 
generated by all elements of even length like the case of the alternating groups. 
We will show this subalgebra has dimension just a half of that of $\mathcal{H}_{n}(q)$. 
Moreover, we will show a presentation of $\frak{A}_{n}(q)$, which is a $q$-analogue of 
that of $\frak{A}_{n}$. Comparing with basic relations of $\frak{A}_{n}$ shown in 
Proposition 3.5, quadratic relations are same as in the case of 
$\frak{A}_{n}$, but cubic relations are slightly different. \par
Now we define the subalgebra $\frak{A}_{n}(q)$ of $\mathcal{H}_{n}(q)$. 
\begin{Def}
For $n \geq 3$ and $q \neq -1$, we define $\frak{A}_{n}(q)$ the subalgebra over $\mathbb{C}$
of $\mathcal{H}_{n}(q)$ 
generated by $f_1f_{i+1}$ with $1 \leq i \leq n-1$.
\end{Def}
To be consistent with the alternating group $\frak{A}_{2}$, we set $\frak{A}_{2}(q)$ the 
subalgebra of $\mathcal{H}_{2}(q)$ generated by the unit element.
\begin{Prop}
$\frak{A}_{n}(q)$ consists of all elements whose terms have even length 
in the expression by $f_i$ with $1 \leq i \leq n-1$.
\end{Prop}
\begin{proof}
From Proposition 4.5, the multiplication of two monomial elements of
even length is expressed by elements whose terms have even length. 
Hence it is enough to see that $f_jf_k$(for $1 \leq j,k \leq n-1$ $j\neq k$) 
are generated by $f_1f_i$'s. 
If $j=1$ and $k>1$ or $j>2$ and $k=1$, then $f_jf_k$ is the generator itself. 
For $j > 2$ and $k > 2$, $f_1f_jf_1f_k=f_jf_k$. 
For $j=2$, we obtain from Proposition 4.5 that 
\[
f_2f_k = (f_1f_2)^2f_1f_k-\Big{(}\dfrac{q-1}{q+1}\Big{)}^2(f_1f_k-f_1f_2f_1f_k)
\]
Thus, we complete the proof.
\end{proof}
\begin{Theo}
As a vector space over the field $\mathbb{C}$, $\frak{A}_{n}(q)$ 
has dimension $n!/2$.
\[ \dim_{\mathbb{C}}\frak{A}_{n}(q)=\dfrac{n!}{2} \]
\end{Theo}
\begin{proof}
Consider all of elements which have even length in the set,
\begin{equation*}
\begin{split}
\{{U'}_{2,d_2}{U'}_{3,d_3}\hdots {U'}_{i,d_i} \hdots {U'}_{n,d_n}\mid & 
(d_2,d_3,\hdots,d_n) \\ 
& \in [0,1]\times[0,2]\times\hdots\times[0,n-1]\} 
\end{split}
\end{equation*}
Then they are linearly independent over $\mathbb{C}$ in $\frak{A}_{n}(q)$. 
Thus it is sufficient to prove that any element in $\frak{A}_{n}(q)$ 
is expressed as a linear combination of the form,
\begin{equation*}
\begin{split}
\{{U'}_{2,d_2}{U'}_{3,d_3} \hdots {U'}_{i,d_i} & \hdots {U'}_{n,d_n}\mid 
(d_2,d_3,\hdots,d_n) \\ 
& \in [0,1]\times[0,2]\times\hdots\times[0,n-1] \quad \text{length is even}
\}
\end{split}
\end{equation*}
We will prove this by induction on 
the length of elements in $\frak{A}_{n}(q)$. 
Consider the element $f_{i_1}f_{i_2}\hdots f_{i_{2k}}$ in $\frak{A}_{n}(q)$.
In the exchange and deletion 
process of expressions of words using relations of Proposition 4.5, 
leading terms of right hand sides of Proposition 4.5 
bring us same calculations as in the case of $\frak{S}_{n}$, 
and remaining terms 
calculations for the products of at most $2(k-1)$ $f_i$'s. 
By the induction assumption, 
remaining terms can be expressed with a linear combination of desired form, 
hence $f_{i_1}f_{i_2}\hdots f_{i_{2k}}$ can be also.
\end{proof}
Next we will give a presentation of $\frak{A}_{n}(q)$.
\begin{Prop}
Let $y_i=f_1f_{i+1}$ for $i=1,2,\hdots ,n-2$. Then $y_i$'s satisfy the 
following relations.
\begin{enumerate}
\rm
\item
\it
$y_1^3=-\Big{(}\dfrac{q-1}{q+1}\Big{)}^2(y_1^2-y_1)+1$
\rm
\item
\it
$y_i^2=1 \qquad \text{for $i>1$}$
\rm
\item
\it
$(y_{i-1}y_i)^3=-\Big{(}\dfrac{q-1}{q+1}\Big{)}^2\Big{\{}(y_{i-1}y_i)^2-y_{i-1}y_i\Big{\}}+1 
\qquad \text{for $i=2,3,\hdots ,n-2$}$
\rm
\item
\it
$(y_iy_j)^2=1 \qquad \text{whenever $|i-j|>1$}$
\end{enumerate}
\end{Prop}
\begin{proof}
Relations 2 and 4 are obvious. For relation 1,
\begin{equation*}
\begin{split}
(y_1)^3 &= f_1f_2f_1f_2f_1f_2 \\
	&= \Big{\{}f_2f_1+\Big{(}\dfrac{q-1}{q+1}\Big{)}^2(1-f_1f_2)\Big{\}}f_1f_2 \\
	&= 1+\Big{(}\dfrac{q-1}{q+1}\Big{)}^2\Big{\{}f_1f_2-(f_1f_2)^2\Big{\}} \\
	&= -\Big{(}\dfrac{q-1}{q+1}\Big{)}^2(y_1^2-y_1)+1
\end{split}
\end{equation*}
For relation 3,
\begin{equation*}
\begin{split}
(y_{i-1}y_i)^3 &= f_1f_if_1f_{i+1}f_1f_if_1f_{i+1}f_1f_if_1f_{i+1} \\
	&= f_1f_if_{i+1}f_if_{i+1}f_if_1f_{i+1} \\
	&= \Big{\{}f_1f_{i+1}f_if_{i+1}
		+\Big{(}\dfrac{q-1}{q+1}\Big{)}^2(f_1f_{i+1}-f_1f_i)\Big{\}}f_{i+1}f_if_1f_{i+1} \\
	&= 1+\Big{(}\dfrac{q-1}{q+1}\Big{)}^2(f_1f_if_1f_{i+1}-f_1f_if_{i+1}f_if_1f_{i+1}) \\
	&= 1+\Big{(}\dfrac{q-1}{q+1}\Big{)}^2\Big{\{}y_{i-1}y_i-(y_{i-1}y_i)^2\Big{\}}
\end{split}
\end{equation*}
Thus we complete the proof.
\end{proof}
\begin{Theo}
The generators $y_i$'s and relations given in Proposition 5.4 define a presentation 
of $\frak{A}_{n}(q)$.
\end{Theo}
We will prove this theorem at the end of this section. \par
Assume that $q \neq -1$. 
Let $\Tilde{\frak{A}}_{n}(q)$($n\geq 3$) be an algebra over $\mathbb{C}$ generated by elements 
$a_1,a_2,\hdots ,a_{n-2}$ with defining relations,
\begin{enumerate}
\item
$a_1^3=-\Big{(}\dfrac{q-1}{q+1}\Big{)}^2(a_1^2-a_1)+1$
\item
$a_i^2=1 \qquad \text{for $i>1$}$
\item
$(a_{i-1}a_i)^3=-\Big{(}\dfrac{q-1}{q+1}\Big{)}^2\Big{\{}(a_{i-1}a_i)^2-a_{i-1}a_i\Big{\}}+1 
\qquad \text{for $i=2,3,\hdots ,n-2$}$
\item
$(a_ia_j)^2=1 \qquad \text{whenever $|i-j|>1$}$
\end{enumerate}
For $n=2$, we set $\Tilde{\frak{A}}_{2}(q)$ an algebra over $\mathbb{C}$ generated by only the unit 
element. \par
To prove the Theorem 5.5, we shall begin with the following Lemma.
\begin{Lem}
The following relations hold for $\Tilde{\frak{A}}_{n}(q)$.
\begin{enumerate}
\rm
\item
\it
$a_1^{-1}=a_1^2+\Big{(}\dfrac{q-1}{q+1}\Big{)}^2(a_1-1)$
\rm
\item
\it
$a_i^{-1}=a_i \qquad \text{for $i=2,3,\hdots ,n-2$}$
\rm
\item
\it
$a_{i-1}a_ia_{i-1}=a_ia_{i-1}^{-1}a_i-\Big{(}\dfrac{q-1}{q+1}\Big{)}^2(a_{i-1}-a_i)\qquad 
\text{for $i=2,3,\hdots ,n-2$}$
\rm
\item
\it
$a_ia_{i-1}a_i=a_{i-1}^{-1}a_ia_{i-1}^{-1}-\Big{(}\dfrac{q-1}{q+1}\Big{)}^2(a_i-a_{i-1}^{-1})\qquad 
\text{for $i=2,3,\hdots ,n-2$}$
\rm
\item
\it
$a_ia_1=a_1^{-1}a_i \qquad a_1a_i=a_ia_1^{-1} \qquad \text{for $i>2$}$
\rm
\item
\it
$a_ia_j=a_ja_i \qquad \text{whenever $|i-j|>1$ and $i,j>1$}$
\end{enumerate}
\end{Lem}
\begin{proof}
A straightforward computation.
\end{proof}
\begin{Prop}
For $n>3$, $\Tilde{\frak{A}}_{n}(q)$ is generated as a vector space over $\mathbb{C}$ by the 
monomials with at most one occurrence of $a_{n-2}$.
\end{Prop}
\begin{proof}
For $n=4$, Let $M$ be a monomial in $\Tilde{\frak{A}}_{4}(q)$ with at least two occurrences of 
$a_2$. Displaying two consecutive occurrences of $a_2$ in $M$, we write $M=M_1a_2M_2a_2M_3$, where 
we can assume that $M_2$ is a monomial in $a_1$, that may be $1$ or $a_1$ or $a_1^2$. 
For the first case, $a_21a_2=1$. For the second case, from the relation 4 of Lemma 5.6,
\[
a_2a_1a_2=a_1^{-1}a_2a_1^{-1}-\Big{(}\dfrac{q-1}{q+1}\Big{)}^2(a_2-a_1^{-1})
\]
Hence we get the desired form. 
For the third case, from the relation 1, 3 of Lemma 5.6,
\begin{equation*}
\begin{split}
a_2a_1^2a_2 
&= a_2\Big{\{}a_1^{-1}-\Big{(}\dfrac{q-1}{q+1}\Big{)}^2(a_1-1)\Big{\}}a_2 \\
&= a_1a_2a_1+\Big{(}\dfrac{q-1}{q+1}\Big{)}^2(a_1-a_2) \\
& \quad -\Big{(}\dfrac{q-1}{q+1}\Big{)}^2(a_2a_1a_2-1)
\end{split}
\end{equation*} 
Applying the second case to the term $a_1a_2a_1$, we get the desired form for the third case. 
Thus we complete the proof for $n=4$. \par
Next, let $M$ be a monomial in $\Tilde{\frak{A}}_{n}(q)$ with at least two occurrences of 
$a_{n-2}$. Displaying two consecutive occurrences of $a_2$ in $M$, we write $M=M_1a_{n-2}M_2a_{n-2}M_3$, 
where we can assume that $M_2$ is a monomial in $a_1,a_2 \hdots a_{n-3}$. 
We can assume by induction that $M_2$ contains $a_{n-3}$ at most once. 
If $M_2$ does not contain $a_{n-3}$ at all, then by the relation 2, 5, 6 of Lemma 5.6,
\[
M=M_1M_2'a_{n-3}^2M_3=M_1M_2'M_3
\]
for some $M_2' \in \Tilde{\frak{A}}_{n-1}(q)$. In this case, $a_{n-3}$ is vanished.
If $M_2$ contains $a_{n-3}$ exactly once, we can write $M_2=M_4a_{n-3}M_5$ with 
$M_4,M_5$ in $a_1,a_2, \hdots a_{n-4}$ and then by the relations 5, 6 of Lemma 5.6,
\[
M=M_1a_{n-2}M_4a_{n-3}M_5a_{n-2}M_3=M_1M_4'a_{n-2}a_{n-3}a_{n-2}M_5'M_3
\]
for some $M_4',M_5' \in \Tilde{\frak{A}}_{n-1}(q)$.
By the relation 4 of Lemma 5.6,
\[
M = M_1M_4'\Big{\{}a_{n-3}^{-1}a_{n-2}a_{n-3}^{-1}-
\Big{(}\dfrac{q-1}{q+1}\Big{)}^2(a_{n-2}-a_{n-3}^{-1})\Big{\}}M_5'M_3
\]
reducing the number of occurrence of $a_{n-2}$.
\end{proof}
Consider the following sets of monomials.
\begin{equation*}
\begin{split}
S_1 &= \{1,a_1,a_1^2\} \\
S_2 &= \{1,a_2,a_2a_1,a_2a_1^2\} \\
\vdots \\
S_i &= \{1,a_i,a_ia_{i-1},\hdots ,a_ia_{i-1}\hdots a_2a_1,a_ia_{i-1}\hdots a_2a_1^2\} \\
\vdots \\
S_{n-2} &= \{1,a_{n-2},a_{n-2}a_{n-3},\hdots ,a_{n-2}a_{n-3}\hdots a_2a_1,a_{n-2}a_{n-3}\hdots a_2a_1^2\}
\end{split}
\end{equation*}
We shall say that $M_0=U_1U_2\hdots U_{n-2}$ is 
a monomial in normal form in $\Tilde{\frak{A}}_{n}(q)$, if $U_i \in S_i$ for $i=1,2,\hdots n-2$. 
There are $n!/2$ of monomials in normal form in $\Tilde{\frak{A}}_{n}(q)$. 
\begin{Prop}
The monomials in normal form in $\Tilde{\frak{A}}_{n}(q)$ generate $\Tilde{\frak{A}}_{n}(q)$ 
as a vector space over $\mathbb{C}$. In particular, 
$\dim_{\mathbb{C}}\Tilde{\frak{A}}_{n}(q) \leq n!/2$ 
\end{Prop}
\begin{proof}
The proof will be by induction. For $n=3$, it is obvious. Let $M$ be an element in 
$\Tilde{\frak{A}}_{n}(q)$. If $M$ contains no $a_{n-2}$, then $M \in 
\Tilde{\frak{A}}_{n-1}(q)$. Hence by the induction assumption, $M$ is expressed by 
a linear combination of monomials in normal form in $\Tilde{\frak{A}}_{n}(q)$. 
If $M$ contains $a_{n-2}$, then by Proposition 5.7, we can write $M=M_1a_{n-2}M_2$ where 
$M_1,M_2 \in \Tilde{\frak{A}}_{n-1}(q)$. 
By induction, $M_2$ is a linear combination of the form $U_1U_2 \hdots U_{n-3}$ with 
$U_i \in S_i$ for $i=1,2, \hdots ,n-3$. By the relations 5, 6 of Lemma 5.6, we have
\[
M_1a_{n-2}U_1U_2 \hdots U_{n-3}=M_1'a_{n-2}U_{n-3}
\]
for some $M_1'$ in $\Tilde{\frak{A}}_{n-1}(q)$. By induction again, $M_1'$ is a linear 
combination of monomials of the form $U_1U_2 \hdots U_{n-3}$ with $U_i \in S_i$. 
Thus $M$ is a linear combination of monomials $U_1U_2 \hdots U_{n-2}$ as desired.
\end{proof}
\renewcommand{\proofname}{Proof of Theorem 5.5}
\begin{proof}
Consider the map $\varphi :\Tilde{\frak{A}}_{n}(q) \longrightarrow \frak{A}_{n}(q)$ 
such that $\varphi (a_i)=y_i$. $\varphi$ defines a homomorphism of algebras. 
We immediately observe that this map is surjective. Hence it is enough to see that 
the dimension of $\Tilde{\frak{A}}_{n}(q)$ is no more than the dimension of $\frak{A}_{n}(q)$, 
but this was already shown in Proposition 5.8. 
\end{proof}
\renewcommand{\proofname}{Proof}

\section{Representations of $\frak{A}_{n}(q)$}\par
In this section, we analyze representations of $\frak{A}_{n}(q)$ by restricting the 
representations we know from $\mathcal{H}_{n}(q)$. 
There is a close relationship between representations of $\frak{S}_{n}$ and 
those of $\frak{A}_{n}$ already shown in section 2. 
We will see this phenomenon in the case of $\mathcal{H}_{n}(q)$ and $\frak{A}_{n}(q)$.
For irreducible representations $\pi_\lambda$ of $\mathcal{H}_{n}(q)$, 
we denote by $\Tilde{\pi}_\lambda$ the restriction of $\pi_\lambda$ to 
$\frak{A}_{n}(q)$. For $g \in \frak{A}_{n}(q)$, 
\[ \Tilde{\pi}_\lambda(g)v_T=\sum_{T'\in \operatorname{STab}(\lambda)}
g^\lambda_{T,T'}v_{T'} \]
where $g^\lambda_{T,T'}$'s are complex numbers.
For $T,T'\in \operatorname{STab}(\lambda)$, there exists $\sigma \in 
\frak{S}_{n}$  such that $\sigma \cdot T = T'$.
We define the distance $d(T,T')$ from $T$ to $T'$ by $l(\sigma)$.
\begin{Lem}
For $g \in \frak{A}_{n}(q)$, 
\[ g^{^t\lambda}_{ ^tT, ^tT'}=g^\lambda_{T,T'} \]
\end{Lem}
\begin{proof}
It is enough to prove for the case $g \in \{y_i \mid 1 \leq i \leq n-2\}$. 
The element $f_i$ acts on $v_T$ in three ways according to where '$i$' and '$i+1$'
appear in 
the standard tableau $T$. Hence, for the action of $y_i=f_1f_{i+1}$, we must consider 
places of $1$, $2$, $i+1$ and $i+2$. We claim that only two cases are possible for 
arrangements of $1$ and $2$: they appear in the same row, or appear in the same column. 
If $i>1$, then there is no interaction 
between the action of $f_1$ and $f_{i+1}$. But if $i=1$, then actions of $f_1$ and $f_2$ 
interact. Hence the action of $y_1$ is more complicate. So, we show equalities of 
matrix elements as the case may be.\par \par
The case $i>1$. \par
\begin{enumerate}
\item
If $1$ and $2$ appear in the same row(resp. column) of $T$, and 
$i+1$ and $i+2$ also appear in the same row(resp. column) of $T$, then 
$1$ and $2$ appear in the same column(resp. row) of $^tT$ and 
$i+1$ and $i+2$ also appear in the same column(resp. row) of $^tT$. 
In this case, we easily get 
$\Tilde{\pi}_\lambda(y_i)v_T=v_T$ and 
$\Tilde{\pi}_{^t\lambda}(y_i)v_{^tT}=v_{^tT}$.
\item
If $1$ and $2$ appear in the same row(resp. column) of $T$, and 
$i+1$ and $i+2$ appear in the same column(resp. row) of $T$, then 
$1$ and $2$ appear in the same column(resp. row) of $^tT$ and 
$i+1$ and $i+2$ appear in the same row(resp. column) of $^tT$. 
In this case, we easily get 
$\Tilde{\pi}_\lambda(y_i)v_T=-v_T$, and 
$\Tilde{\pi}_{^t\lambda}(y_i)v_{^tT}=-v_{^tT}$.
\item
If $1$ and $2$ appear in the same row(resp. column) of $T$, and 
$i+1$ and $i+2$ appear in neither the same row nor the same column 
of $T$, then $1$ and $2$ appear in the same column(resp. row) of $^tT$. 
Using the notation of the end of section 4, set 
\[ \Tilde{\pi}_\lambda(f_{i+1})=
\begin{bmatrix}
	A & B \\
	B & -A \\
\end{bmatrix} \qquad \text{for the basis $v_T,v_{s_{i+1}\cdot{T}}$}
\]
and we get
\begin{equation*}
\begin{split}
\Tilde{\pi}_\lambda(y_i)v_T &= Av_T+Bv_{s_{i+1}\cdot T} 
\text{(resp. $-Av_{T}-Bv_{s_{i+1}\cdot{T}}$)} \\
\Tilde{\pi}_{^t\lambda}(y_i)v_{^tT} &= Av_{^tT}+Bv_{s_{i+1}\cdot{^tT}} 
\text{(resp. $-Av_{^tT}-Bv_{s_{i+1}\cdot{^tT}}$)}
\end{split}
\end{equation*}
\end{enumerate}\par
The case $i=1$ \par
By exchanging $T$ and $^tT$, we observe that it is enough to check only in two cases: 
$1$, $2$ and $3$ appear in the same row, 
$1$ and $2$ appear in the same row and $1$ and $3$ appear in the same column.
\begin{enumerate}
\item
If $1$, $2$ and $3$ appear in the same row of $T$, 
then $1$, $2$ and $3$ appear in the same column of $^tT$. 
In this case, we easily get 
$\Tilde{\pi}_\lambda(y_1)v_T=v_T$ and 
$\Tilde{\pi}_{^t\lambda}(y_1)v_{^tT}=v_{^tT}$.
\item
$1$ and $2$ appear in the same row and $1$ and $3$ 
appear in the same column of $T$, then $1$ and $2$ appear in the same column and $1$ and $3$ 
appear in the same row of $^tT$.
We observe that $f_1$ acts on $v_{T}$, 
$v_{s_2\cdot{T}}$, $v_{^tT}$ 
and $v_{s_2{^tT}}$ as scalar multiplication. 
Using the notation of the end of section 4, set 
\begin{equation*}
\begin{split}
\Tilde{\pi}_\lambda(f_2) &= 
\begin{bmatrix}
	A & B \\
	B & -A \\
\end{bmatrix} \qquad \text{for the basis 
$v_{T},v_{s_2\cdot{T}}$} \\
\end{split}
\end{equation*}
and we get
\begin{equation*}
\begin{split}
\Tilde{\pi}_\lambda(y_1)v_T &= Av_T - Bv_{s_2\cdot T} \\
\Tilde{\pi}_{^t\lambda}(y_1)v_{^tT} &= Av_{^tT} - Bv_{s_2\cdot ^tT} 
\end{split}
\end{equation*}
\end{enumerate}
\end{proof}
\begin{Prop}
For $\lambda \in \Lambda_n$, 
\[ \Tilde{\pi}_\lambda \cong \Tilde{\pi}_{^t\lambda} \]
\end{Prop}
\begin{proof}
From Lemma 6.1, representation matrices of $\Tilde{\pi}_\lambda$ and $\Tilde{\pi}_{^t\lambda}$ 
for generators of $\frak{A}_{n}(q)$ are coincide. Therefore, the above statement holds.
\end{proof}
For $\lambda \in \Lambda_n$ with $n>1$, we set
\begin{equation*}
\begin{split}
\operatorname{STab}(\lambda)^+ &= \{T \in \operatorname{STab}(\lambda) \mid T(1,1)=1,T(1,2)=2\} \\
\operatorname{STab}(\lambda)^- &= \{T \in \operatorname{STab}(\lambda) \mid T(1,1)=1,T(2,1)=2\} 
\end{split}
\end{equation*}
Then, $\operatorname{STab}(\lambda)$ is a disjoint union of 
$\operatorname{STab}(\lambda)^+$ and 
$\operatorname{STab}(\lambda)^-$. If $\lambda$ is self-conjugate, then 
$T \in \operatorname{STab}(\lambda)^+$ if and only if ${^tT} \in 
\operatorname{STab}(\lambda)^-$, 
therefore, $\operatorname{STab}(\lambda)^+$ corresponds to 
$\operatorname{STab}(\lambda)^-$ one to one and 
onto by mapping, $T \longrightarrow {^tT}$. \par
Let $\lambda$ is self-conjugate. 
Let $\Tilde{V}_\lambda$ be the 
representation space of $\Tilde{\pi}_\lambda$. Now, we will introduce 
two subspaces $\Tilde{V}_\lambda^+$ and $\Tilde{V}_\lambda^-$ of 
$\Tilde{V}_\lambda$ as follows. 
\begin{equation*}
\begin{split}
\Tilde{V}_\lambda^+ &= \bigoplus_{T \in \operatorname{STab}(\lambda)^+}
\mathbb{C}(v_{T'}+v_{^tT}) \\
\Tilde{V}_\lambda^- &= \bigoplus_{T \in \operatorname{STab}(\lambda)^+}
\mathbb{C}(v_{T}-v_{^tT})
\end{split}
\end{equation*}
Then, $\Tilde{V}_\lambda$ is a direct sum of $\Tilde{V}_\lambda^+$ and 
$\Tilde{V}_\lambda^-$ as vector space over $\mathbb{C}$.
\begin{Prop}
Let $\lambda \in \Lambda_n$ be self-conjugate. 
Then $\Tilde{V}_\lambda^+$ and $\Tilde{V}_\lambda^-$ are $\frak{A}_{n}(q)$-
submodules of $\Tilde{V}_\lambda$, and $\Tilde{V}_\lambda$ has a 
$\frak{A}_{n}(q)$-submodule decomposition as follows.
\[ \Tilde{V}_\lambda=\Tilde{V}_\lambda^+ \oplus 
\Tilde{V}_\lambda^- \qquad \text{as $\frak{A}_{n}(q)$-module} \]
\end{Prop}
\begin{proof}
Let $g \in \frak{A}_{n}(q)$. 
We set $T \in \operatorname{STab}(\lambda)^+$. 
From Lemma 6.1, we get following calculation for the element $v_{T}+v_{^tT}$ 
of $\Tilde{V}_\lambda^+$.
\begin{equation*}
\begin{split}
\Tilde{\pi}_\lambda(g)(v_{T}+v_{^tT}) 
	&= \sum_{T' \in \operatorname{STab}(\lambda)}g_{TT'}^{\lambda}(v_{T'}+v_{^tT'}) \\
	&= \sum_{T' \in \operatorname{STab}(\lambda)^+}g_{TT'}^{\lambda}(v_{T'}+v_{^tT'}) \\
	&\quad + \sum_{T' \in \operatorname{STab}(\lambda)^-}g_{TT'}^{\lambda}(v_{T'}+v_{^tT'}) \\
	&= \sum_{T' \in \operatorname{STab}(\lambda)^+}g_{TT'}^{\lambda}(v_{T'}+v_{^tT'}) \\
	&\quad + \sum_{T' \in \operatorname{STab}(\lambda)^+}g_{T{^tT'}}^{\lambda}(v_{T'}+v_{^tT'}) \\
	&= \sum_{T' \in \operatorname{STab}(\lambda)^+}
		(g_{TT'}^{\lambda}+g_{T{^tT'}}^{\lambda})(v_{T'}+v_{^tT'}) \\
\end{split}
\end{equation*}
Thus, $\Tilde{\pi}_\lambda(g)(v_{T'}+v_{^tT})$ is in $\Tilde{V}_\lambda^+$. \par
Next, we get following calculation for the element $v_{T'}-v_{^tT'}$ 
of $\Tilde{V}_\lambda^-$.
\begin{equation*}
\begin{split}
\Tilde{\pi}_\lambda(g)(v_{T}-v_{^tT}) 
	&= \sum_{T' \in \operatorname{STab}(\lambda)}g_{TT'}^{\lambda}(v_{T'}-v_{^tT'}) \\
	&= \sum_{T' \in \operatorname{STab}(\lambda)^+}g_{TT'}^{\lambda}(v_{T'}-v_{^tT'}) \\
	&\quad + \sum_{T' \in \operatorname{STab}(\lambda)^-}g_{TT'}^{\lambda}(v_{T'}-v_{^tT'}) \\
	&= \sum_{T' \in \operatorname{STab}(\lambda)^+}g_{TT'}^{\lambda}(v_{T'}-v_{^tT'}) \\
	&\quad - \sum_{T' \in \operatorname{STab}(\lambda)^+}g_{T'{^tT'}}^{\lambda}(v_{T'}-v_{^tT'}) \\
	&= \sum_{T' \in \operatorname{STab}(\lambda)^+}
		(g_{TT'}^{\lambda}-g_{T{^tT'}}^{\lambda})(v_{T'}-v_{^tT'}) \\
\end{split}
\end{equation*}
Thus, $\Tilde{\pi}_\lambda(g)(v_{T}-v_{^tT})$ is in $\Tilde{V}_\lambda^-$. \\
\end{proof}
We denote by $\Tilde{\pi}_{\lambda}^+$ the representation of $\frak{A}_{n}(q)$ corresponding to
$\Tilde{V}_{\lambda}^+$ and $\Tilde{\pi}_{\lambda}^-$ corresponding to $\Tilde{V}_{\lambda}^-$.\par
Next, we will show the irreducibilities of these representations and semisimplicity 
of $\frak{A}_{n}(q)$. 
Let $\lambda \in \Lambda_n$ be non self-conjugate. Then $\Tilde{V}_\lambda$ and 
$\Tilde{V}_{^t\lambda}$ are isomorphic as $\frak{A}_{n}(q)$-modules and 
the direct sum $\Tilde{V}_\lambda \oplus \Tilde{V}_{^t\lambda}$ has 
another decomposition as $\mathbb{C}$ vector spaces,
\[ \Tilde{V}_\lambda \oplus \Tilde{V}_{^t\lambda} 
	= \Tilde{V}_\lambda^+ \oplus \Tilde{V}_\lambda^- \]
where $\Tilde{V}_\lambda^+$ and $\Tilde{V}_\lambda^-$ defined as follows.
\begin{equation*}
\begin{split}
\Tilde{V}_\lambda^+ &= \bigoplus_{T \in \operatorname{STab}(\lambda)}
			\mathbb{C}(v_{T}+v_{^tT}) \\
\Tilde{V}_\lambda^- &= \bigoplus_{T \in \operatorname{STab}(\lambda)}
			\mathbb{C}(v_{T}-v_{^tT}) 
\end{split}
\end{equation*}
\begin{Prop}
If $\lambda$ is non self-conjugate, then the above decomposition of 
$\Tilde{V}_\lambda \oplus \Tilde{V}_{^t\lambda}$,
\[ \Tilde{V}_\lambda \oplus \Tilde{V}_{^t\lambda} 
	= \Tilde{V}_\lambda^+ \oplus \Tilde{V}_\lambda^- \]
is a submodule decomposition of $\frak{A}_{n}(q)$.
\end{Prop}
\begin{proof}
For $v_{T}+v_{^tT} \in \Tilde{V}_\lambda^+ \subset 
\Tilde{V}_\lambda \oplus \Tilde{V}_{^t\lambda}$, we consider the direct sum 
of irreducible representations $\Tilde{\pi}_\lambda$ and $\Tilde{\pi}_{^t\lambda}$.
For $g \in \frak{A}_{n}(q)$, we have from Lemma 6.1 the following.
\begin{equation*}
\begin{split}
(\Tilde{\pi}_\lambda\oplus\Tilde{\pi}_{^t\lambda})(g)(v_{T}+v_{^tT}) 
	&= \sum_{T' \in \operatorname{STab}(\lambda)}g_{TT'}^{\lambda}(v_{T'}+v_{^tT'}) \\
\end{split}
\end{equation*}
Thus, ($\Tilde{\pi}_{\lambda}\oplus\Tilde{\pi}_{^t\lambda})(g)
(v_{T}+v_{^tT})$ is in $\Tilde{V}_\lambda^+$. \par
The Similar calculation is valid for $\Tilde{V}_\lambda^-$, hence 
we omit the calculation for $\Tilde{V}_\lambda^-$.
\end{proof}
\begin{Theo}
Let $q$ be a complex number such that $q\neq 0$ and $q$ is not a 
$k$-th root of unity with $1\leq k \leq n$. 
Let $\lambda_1$,$^t\lambda_1$,$\lambda_2$,$^t\lambda_2$,$\hdots$,
$\lambda_p$,$^t\lambda_p$ be non self-conjugate 
Young diagrams and $\lambda_{p+1}$,$\lambda_{p+2}$,$\hdots$,$\lambda_{p+q}$ 
be self-conjugate Young diagrams within $\Lambda_n$. 
Then representations 
$\Tilde{\pi}_{\lambda_1}$,$\Tilde{\pi}_{\lambda_2}$,
$\hdots,\Tilde{\pi}_{\lambda_p}$,$\Tilde{\pi}_{\lambda_{p+1}}^+$,
$\Tilde{\pi}_{\lambda_{p+1}}^-$,$\hdots,\Tilde{\pi}_{\lambda_{p+q}}^+$,
$\Tilde{\pi}_{\lambda_{p+q}}^-$ are irreducible and not equivalent 
each other. 
These representations consist of all equivalent classes of 
irreducible representations of $\frak{A}_{n}(q)$. Hence, 
$\frak{A}_{n}(q)$ is semisimple.
\end{Theo}
\begin{proof}
At first, we will show the semisimplicity of $\frak{A}_{n}(q)$ 
under the assumptions of irreducibilities and mutual inequalities 
of these representations. We consider the map 
\begin{equation*}
\begin{split}
\Tilde{\pi}_n : x\in \frak{A}_{n}(q) \longrightarrow 
\Tilde{\pi}_{\lambda_1}(x)\oplus \hdots \oplus \Tilde{\pi}_{\lambda_p}(x) 
\oplus \Tilde{\pi}_{\lambda_{p+1}}^+(x) & \oplus \Tilde{\pi}_{\lambda_{p+1}}^-(x) 
\oplus \hdots \\ 
& \hdots \oplus \Tilde{\pi}_{\lambda_{p+q}}^+(x) \oplus \Tilde{\pi}_{\lambda_{p+q}}^-(x)
\end{split}
\end{equation*} 
Then, by theorems of Burnside and Frobenius-Schur, 
$\frak{A}_{n}(q)$ has a quotient $\Tilde{\pi}_n(\frak{A}_{n}(q))$ 
isomorphic to the semisimple algebra 
$\oplus \operatorname{End}_{\mathbb{C}}V$ where $V$ runs over irreducible 
representation spaces listed in the statement of theorem.\par
This semisimple algebra has dimension $n!/2$ and we already show that 
$\frak{A}_{n}(q)$ has dimension $n!/2$, thus $\frak{A}_{n}(q)$ 
is isomorphic to $\oplus \operatorname{End}_{\mathbb{C}}V$ and semisimple.\par
Next, we will show the irreducibilities and mutual inequalities 
of these representations by induction.\par
For $n=2$, it is obvious. Indeed, $\frak{A}_{2}(q)$ is generated by the unit element, 
and both $\Tilde{\pi}_{(2)}$ and $\Tilde{\pi}_{(1^2)}$ are identity maps. \par
For $n=3$, $\Tilde{\pi}_{(3)}$ and $\Tilde{\pi}_{(1^3)}$ are identity map and 
\[ \Tilde{V}_{(2,1)}=\Tilde{V}_{(2,1)}^+ \oplus \Tilde{V}_{(2,1)}^- \]
where, 
\[ \Tilde{V}_{(2,1)}^+=\mathbb{C}(v_{T}+v_{^tT}) \]
\[ \Tilde{V}_{(2,1)}^-=\mathbb{C}(v_{T}-v_{^tT}) \]
with standard tableau $T$ as follows.
\[
T(1,1)=1 \qquad T(1,2)=2 \qquad T(2,1)=3
\]
%\begin{figure}[h]
%\begin{center}
%\setlength{\unitlength}{1mm}
%\begin{picture}(20,10)
%\put(0,5){$T=$}
%\put(10,0){\framebox(5,5){$3$}}
%\put(10,5){\framebox(5,5){$1$}}
%\put(15,5){\framebox(5,5){$2$}}
%\end{picture}
%\end{center}
%\end{figure} \par
We get the following from the proof of Lemma 6.1,Proposition 6.3,
\begin{equation*}
\begin{split}
\Tilde{\pi}_{(2,1)}^\pm (y_1)(v_{T}\pm v_{^tT}) &=
	\dfrac{1+q^2 \mp 2\sqrt{q(1+q+q^2)}}{(1+q)^2}(v_{T}\pm v_{^tT}) \\
\Tilde{\pi}_{(2,1)}^\pm (y_1^2)(v_{T}\pm v_{^tT}) &=
	\Big{\{}\dfrac{1+q^2 \mp 2\sqrt{q(1+q+q^2)}}{(1+q)^2}\Big{\}}^2(v_{T}\pm v_{^tT}) \\
\end{split}
\end{equation*}
Since every representation has degree $1$, these representations are 
irreducible. Moreover, since $q$ is neither $0$ nor $k$-th root of unity with 
$1 \leq k \leq 3$, we immediately have $\sqrt{q(1+q+q^2)}$ is nonzero. 
Hence they are mutually inequivalent.\par
Let $n>3$. 
By induction assumption, $\frak{A}_{n-1}(q)$ is a semisimple algebra with 
central primitive idempotents $z_{{\lambda '}_1}$,$z_{{\lambda '}_2}$,$\hdots$,
$z_{{\lambda '}_{p'}}$,$z_{{\lambda '}_{p'+1}}^+$,$z_{{\lambda '}_{p'+1}}^-$,$\hdots$,
$z_{{\lambda '}_{p'+q'}}^+$,$z_{{\lambda '}_{p'+q'}}^-$ with 
${\lambda '}_1$,$\hdots$,${\lambda '}_{p'+q'}$$\in \Lambda_{n-1}$ and 
${\lambda '}_1$,$\hdots$,${\lambda '}_{p'}$ are non self-conjugate and 
${\lambda '}_{p'+1}$,$\hdots$, 
${\lambda '}_{p'+q'}$ are self-conjugate. \par
Let $\lambda \in \Lambda_n$ be non self-conjugate. Then there is no diagram 
$\lambda '\in \Lambda_{n-1}$ such that $\lambda '$ is non self-conjugate and 
$\lambda '<\lambda$ and 
${^t\lambda '}< \lambda$. Indeed, for such diagram, we immediately obtain,
${^t\lambda '}<{^t\lambda}$ and ${\lambda '}<{^t\lambda}$. Hence, 
${^t\lambda = \lambda}$ but this is contradiction.
When we restrict $\Tilde{V}_\lambda$ to $\frak{A}_{n-1}(q)$-module, 
we can write(for the reason, see [15]), 
\begin{equation*}
\Tilde{V}_{\lambda} = \Tilde{V}_{{\lambda '}_{i_1}} \oplus \Tilde{V}_{{\lambda '}_{i_2}} 
\oplus \hdots \oplus \Tilde{V}_{{\lambda '}_{i_r}} \oplus \Tilde{V}_{{\lambda '}_{i_{r+1}}} 
\oplus \hdots \oplus \Tilde{V}_{{\lambda '}_{i_{r+s}}}
\end{equation*}
where ${\lambda '}_{i_j}$'s are all of elements in $\Lambda_n$ such that 
${\lambda '}_{i_j} < \lambda$. 
We suppose that ${\lambda '}_{i_1},{\lambda '}_{i_2},\hdots ,{\lambda '}_{i_r}$ are 
non self-conjugate and 
${\lambda '}_{i_{r+1}},\hdots ,{\lambda '}_{i_{r+s}}$ are self-conjugate. We observe that 
$s$ is at most $1$ because it is impossible to remove one box and add another box 
with keeping the self-conjugacy. We can write by induction,
\begin{equation*}
\Tilde{V}_{\lambda} = \Tilde{V}_{{\lambda '}_{i_1}} \oplus \Tilde{V}_{{\lambda '}_{i_2}} 
\oplus \hdots \oplus \Tilde{V}_{{\lambda '}_{i_r}}
\end{equation*}
or
\begin{equation*}
\Tilde{V}_{\lambda} = \Tilde{V}_{{\lambda '}_{i_1}} \oplus \Tilde{V}_{{\lambda '}_{i_2}} 
\oplus \hdots \oplus \Tilde{V}_{{\lambda '}_{i_r}} \oplus \Tilde{V}_{{\lambda '}_{i_{r+1}}}^+ 
\oplus \Tilde{V}_{{\lambda '}_{i_{r+1}}}^-
\end{equation*}
with each subspace is irreducible $\frak{A}_{n-1}(q)$-module and inequivalent each other. 
For ${\lambda}',\Tilde{\lambda}' \in \Lambda_{n-1}$ such that ${\lambda}'<\lambda$ and 
$\Tilde{\lambda}'<\lambda$, there is exactly one ${\lambda}'' \in \Lambda_{n-2}$ such that 
${\lambda}'' <{\lambda}'$ and ${\lambda}'' <\Tilde{\lambda}'$. Let $T \in 
\operatorname{STab}(\lambda)$ be such that the tableau obtained from $T$ by removing 
$n$-th box is shape ${\lambda}'$ and the tableau obtained from $T$ by removing 
$n$-th box and $n-1$-th box is shape ${\lambda}''$. Similarly, let $\Tilde{T} \in 
\operatorname{STab}(\lambda)$ be such that the tableau obtained from $\Tilde{T}$ by removing 
$n$-th box is shape $\Tilde{\lambda}'$ and the tableau obtained from $\Tilde{T}$ 
by removing $n$-th box and $n-1$-th box is shape ${\lambda}''$. \par
If ${\lambda}'$ is non self-conjugate, then $v_T \in \Tilde{V}_{{\lambda}'}$, and we get,
\[ \Tilde{\pi}_{\lambda}(y_{n-2})v_T=Av_T+Bv_{\Tilde{T}} \]
Since $q\neq 0$ and $q$ is not a $k$-th root of unity with $1\leq k \leq n$, 
$B$ is well-defined and nonzero. 
If $\Tilde{\lambda}'$ is non self-conjugate, then
\[ \Tilde{\pi}_{\lambda}(z_{\Tilde{\lambda}'})\Tilde{\pi}_{\lambda}(y_{n-2})v_T=
	Bv_{\Tilde{T}} \in \Tilde{V}_{\Tilde{\lambda}'} \]
If $\Tilde{\lambda}'$ is self-conjugate, then as the submodule decomposition
\[ \Tilde{V}_{\Tilde{\lambda}'}
	=\Tilde{V}_{\Tilde{\lambda}'}^+ \oplus \Tilde{V}_{\Tilde{\lambda}'}^- \]
we can write $v_{\Tilde{T}}$ as follows
\[ v_{\Tilde{T}}=
\dfrac{1}{2}(v_{\Tilde{T}}+v_{\Tilde{\Tilde{T}}})+\dfrac{1}{2}(v_{\Tilde{T}}-v_{\Tilde{\Tilde{T}}}) 
\]
where $\Tilde{\Tilde{T}} \in 
\operatorname{STab}(\lambda)$ be such that the tableau obtained from $\Tilde{\Tilde{T}}$ 
by removing 
$n$-th box is shape $\Tilde{\lambda}'$ and the tableau obtained from $\Tilde{\Tilde{T}}$ 
by removing $n$-th box and $n-1$-th box is shape ${^t{{\lambda}''}}$. 
$v_{\Tilde{T}}+v_{\Tilde{\Tilde{T}}} \in \Tilde{V}_{\Tilde{\lambda}'}^+$ and 
$v_{\Tilde{T}}-v_{\Tilde{\Tilde{T}}} \in \Tilde{V}_{\Tilde{\lambda}'}^-$. Hence
\begin{equation*}
\begin{split}
\Tilde{\pi}_{\lambda}(z_{\Tilde{\lambda}'}^+)\Tilde{\pi}_{\lambda}(y_{n-2})v_T &= 
	\dfrac{B}{2}(v_{\Tilde{T}}+v_{\Tilde{\Tilde{T}}}) 
	\in \Tilde{V}_{\Tilde{\lambda}'}^+ \\
\Tilde{\pi}_{\lambda}(z_{\Tilde{\lambda}'}^-)\Tilde{\pi}_{\lambda}(y_{n-2})v_T &= 
	\dfrac{B}{2}(v_{\Tilde{T}}-v_{\Tilde{\Tilde{T}}}) 
	\in \Tilde{V}_{\Tilde{\lambda}'}^- 
\end{split}
\end{equation*}
If ${\lambda}'$ is self-conjugate, then other Young diagrams $\Tilde{\lambda}' \in 
\Lambda_{n-1}$ with 
$\Tilde{\lambda}'<\lambda$ are non self-conjugate. 
We set $v_T+v_{\Hat{T}} \in \Tilde{V}_{{\lambda}'}^+$ 
and $v_T-v_{\Hat{T}} \in \Tilde{V}_{{\lambda}'}^-$ where $\Hat{T} \in 
\operatorname{STab}(\lambda)$ be such that the tableau obtained from $\Hat{T}$ 
by removing 
$n$-th box is shape ${\lambda}'$ and the tableau obtained from $\Hat{T}$ 
by removing $n$-th box and $n-1$-th box is shape ${^t{{\lambda}''}}$. 
Then we get,
\[ \Tilde{\pi}_{\lambda}(y_{n-2})(v_T+v_{\Hat{T}})
	=Av_T+Bv_{\Tilde{T}}
	+\Tilde{\pi}_{\lambda}(y_{n-2})(v_{\Hat{T}}) \]
and elements in $\Tilde{V}_{\Tilde{\lambda}'}$ do not appear in 
$\Tilde{\pi}_{\lambda}(y_{n-2})(v_{\Hat{T}})$. Hence,
\[ \Tilde{\pi}_{\lambda}(z_{\Tilde{\lambda}'})\Tilde{\pi}_{\lambda}(y_{n-2})
	(v_T+v_{\Hat{T}})=Bv_{\Tilde{T}} 
	\in \Tilde{V}_{\Tilde{\lambda}'} \]
is nonzero. \par
The same discussion is valid for $v_T-v_{\Hat{T}} \in \Tilde{V}_{{\lambda}'}^-$ and 
we get,
\[ \Tilde{\pi}_{\lambda}(z_{\Tilde{\lambda}'})\Tilde{\pi}_{\lambda}(y_{n-2})
	(v_T-v_{\Hat{T}})=Bv_{\Tilde{T}} 
	\in \Tilde{V}_{\Tilde{\lambda}'} \]
Because $\Tilde{\lambda}' \in \Lambda_{n-1}$ is arbitrary with $\Tilde{\lambda}'<\lambda$, 
every $\frak{A}_{n}(q)$-submodule $W$ of $\Tilde{V}_{\lambda}$ includes 
whole $\Tilde{V}_{\lambda}$, therefore $\Tilde{V}_{\lambda}$ is irreducible. \par
Let $\lambda \in \Lambda_n$ be self-conjugate. 
When we restrict $\Tilde{V}_\lambda$ to $\frak{A}_{n-1}(q)$-module, 
we can write(see [15] again), 
\begin{equation*}
\Tilde{V}_{\lambda} = \Tilde{V}_{{\lambda '}_{i_1}} \oplus \Tilde{V}_{^t{\lambda '}_{i_1}} 
\oplus \hdots \oplus \Tilde{V}_{{\lambda '}_{i_r}} \oplus \Tilde{V}_{^t{\lambda '}_{i_r}} 
\oplus \Tilde{V}_{{\lambda '}_{i_{r+1}}} 
\oplus \hdots \oplus \Tilde{V}_{{\lambda '}_{i_{r+s}}}
\end{equation*}
where ${\lambda '}_{i_j}$'s are same as the non self-conjugate case. We observe that 
$s$ is at most $1$ because it is impossible to remove one box and add another box 
with keeping the self-conjugacy. 
We set for $j=+,-$,
\[ \Tilde{V}_{\lambda}^j=\Tilde{V}_{{\lambda '}_{i_1}}^j \oplus \hdots \oplus
		\Tilde{V}_{{\lambda '}_{i_r}}^j \oplus \hdots \oplus 
		\Tilde{V}_{{\lambda '}_{i_{r+s}}}^j
\]
where each direct summand is irreducible $\frak{A}_{n-1}(q)$-module with 
$\Tilde{V}_{{\lambda '}_{i_1}}^+ \cong \Tilde{V}_{{\lambda '}_{i_1}}^-$,$\hdots$,
$\Tilde{V}_{{\lambda '}_{i_r}}^+ \cong \Tilde{V}_{{\lambda '}_{i_r}}^-$ and (if exist) 
$\Tilde{V}_{{\lambda '}_{i_{r+1}}}^+$, $\Tilde{V}_{{\lambda '}_{i_{r+1}}}^-$ 
are mutually inequivalent each other by induction.
Hence, from Proposition 6.3, we can write,
\begin{equation*}
\Tilde{V}_{\lambda} = \Tilde{V}_{\lambda}^+ \oplus \Tilde{V}_{\lambda}^- 
\end{equation*}
For ${\lambda}',\Tilde{\lambda}' \in \Lambda_{n-1}$ such that ${\lambda}'<\lambda$ and 
$\Tilde{\lambda}'<\lambda$, there is exactly one ${\lambda}'' \in \Lambda_{n-2}$ such that 
${\lambda}'' <{\lambda}'$ and ${\lambda}'' <\Tilde{\lambda}'$. 
Let $T \in \operatorname{STab}(\lambda)$ be such that the tableau obtained from $T$ by removing 
$n$-th box is shape ${\lambda}'$ and the tableau obtained from $T$ by removing 
$n$-th box and $n-1$-th box is shape ${\lambda}''$. Similarly, let $\Tilde{T} \in 
\operatorname{STab}(\lambda)$ be such that the tableau obtained from $\Tilde{T}$ by removing 
$n$-th box is shape $\Tilde{\lambda}'$ and the tableau obtained from $\Tilde{T}$ 
by removing $n$-th box and $n-1$-th box is shape ${\lambda}''$. \par 
For $v_{T}+v_{^tT} \in \Tilde{V}_{\lambda '}^+$, 
\begin{equation*}
\begin{split}
\Tilde{\pi}_{\lambda}(y_{n-2})(v_{T}+v_{^tT})
&= Av_{T}+Bv_{\Tilde{T}}+Av_{^tT}+Bv_{^t\Tilde{T}} \\
&= A(v_{T}+v_{^tT})+B(v_{\Tilde{T}}+v_{^t\Tilde{T}})
\end{split}
\end{equation*}
where $B(v_{\Tilde{T}}+v_{^t\Tilde{T}})$ is a 
nonzero element in $\Tilde{V}_{\Tilde{\lambda}'}^+$. Hence,
\[ \Tilde{\pi}_{\lambda}(z_{\Tilde{\lambda}'})\Tilde{\pi}_{\lambda}(y_{n-2})
	(v_{T}+v_{^tT})=B(v_{\Tilde{T}}+v_{^t\Tilde{T}}) 
	\in \Tilde{V}_{\Tilde{\lambda}'}^+ 
\]
Because $\Tilde{\lambda}' \in \Lambda_{n-1}$ is arbitrary with $\Tilde{\lambda}'<\lambda$, 
every $\frak{A}_{n}(q)$-submodule $W$ of $\Tilde{V}_{\lambda}^+$ includes 
whole $\Tilde{V}_{\lambda}^+$, hence $\Tilde{V}_{\lambda}^+$ is irreducible.
The same argument is valid for the case in $v_{T}-v_{^tT} \in \Tilde{V}_{\lambda}^-$. 
Thus the irreducibilities of $\Tilde{V}_{\lambda}^+$ and $\Tilde{V}_{\lambda}^-$ were proved. \par
Finally, we show inequivalencies of $\Tilde{V}_{\lambda}$'s as $\frak{A}_{n}(q)$-module. 
If $\lambda$ is self-conjugate, then there is at most one self-conjugate ${\lambda'}$ in 
$\Lambda_{n-1}$ or $\Lambda_{n-2}$ such that ${\lambda'}<{\lambda}$. 
In this case, $\Tilde{V}_{\lambda}^+$ has a direct summand $\Tilde{V}_{\lambda '}^+$ 
and does not have $\Tilde{V}_{\lambda '}^-$ as $\frak{A}_{n-1}(q)$ nor $\frak{A}_{n-2}(q)$-
module, and $\Tilde{V}_{\lambda}^-$ has a direct summand $\Tilde{V}_{\lambda '}^-$ 
and does not have $\Tilde{V}_{\lambda '}^+$ as $\frak{A}_{n-1}(q)$ nor $\frak{A}_{n-2}(q)$-
module  On the other hand, If $\mu$ is non self-conjugate, then $\Tilde{V}_{\mu}$ 
has both $\Tilde{V}_{\Bar\lambda '}^+$ and $\Tilde{V}_{\lambda '}^-$ or has neither 
$\Tilde{V}_{\lambda '}^+$ nor $\Tilde{V}_{\lambda '}^-$. Hence, 
$\Tilde{V}_{\lambda}^+$ and $\Tilde{V}_{\lambda}^-$ and $\Tilde{V}_{\mu}$ are 
mutually inequivalent. \par
Next, we consider the case that different Young diagrams $\lambda$ and 
$\mu$ are both self-conjugate or the case that different Young diagrams 
$\lambda$ and $\mu$ are both non self-conjugate with ${^t\lambda} \neq {\mu}$.
In the former case, the set of all ${\lambda '} \in 
\Lambda_{n-1}$ such that ${\lambda '}<{\lambda}$ do not coincide with 
the set of all ${\mu '} \in \Lambda_{n-1}$ such that ${\mu '}<{\mu}$. 
Hence, $\Tilde{V}_{\lambda}$ and $\Tilde{V}_{\mu}$ 
are inequivalent already as $\frak{A}_{n-1}(q)$-module. 
Since $\Tilde{V}_{\lambda}^+\cong \Tilde{V}_{\lambda}^-$ as $\frak{A}_{n-1}(q)$-module 
and $\Tilde{V}_{\lambda}=\Tilde{V}_{\lambda}^+ \oplus \Tilde{V}_{\lambda}^-$ as 
$\frak{A}_{n-1}(q)$-module, $\Tilde{V}_{\lambda}^\pm$ and $\Tilde{V}_{\mu}^\pm$ are 
mutually inequivalent already as $\frak{A}_{n-1}(q)$-module. 
In the latter case, the set of all ${\lambda '},{^t\lambda '} \in \Lambda_{n-1}$ such that 
${\lambda '}<{\lambda}$ do not coincide with 
the set of all ${\mu '},{^t\mu '} \in \Lambda_{n-1}$ such that ${\mu '}<{\mu}$. Hence, 
$\Tilde{V}_{\lambda}$ and $\Tilde{V}_{\mu}$ are 
mutually inequivalent already as $\frak{A}_{n-1}(q)$-module.
\end{proof}
Finally, we shall show the multiplicity of the irreducible $\mathcal{H}_{n}(q)$-module 
${V}_{\mu}$ in the induced module 
$\operatorname{Ind}_{\frak{A}_{n}(q)}^{\mathcal{H}_{n}(q)}(\Tilde{V})$ where $\Tilde{V}$ 
is an irreducible $\frak{A}_{n}(q)$-module. We use the following Proposition known as 
Nakayama relation.
\begin{Prop}
If $S$ is a subring of a ring $R$, $M$ is a $R$-module and $N$ is a $S$-module then there 
is a natural homomorphism.
\[
\operatorname{Hom}_S(N,\operatorname{Res}_{S}^{R}(M)) \cong 
\operatorname{Hom}_R(\operatorname{Ind}_{S}^{R}(N),M)
\]
\end{Prop}
\begin{proof}
See for example, [1].
\end{proof}
Using the Proposition 6.6 and Schur's Lemma, We have the following.
\begin{Prop}
If $\lambda \in \Lambda_n$ is non self-conjugate then the following holds.
\[
\operatorname{Ind}_{\frak{A}_{n}(q)}^{\mathcal{H}_{n}(q)}(\Tilde{V}_{\lambda}) 
\cong V_{\lambda} \oplus V_{^t\lambda}
\]
If $\lambda \in \Lambda_n$ is self-conjugate then the following holds.
\[
\operatorname{Ind}_{\frak{A}_{n}(q)}^{\mathcal{H}_{n}(q)}(\Tilde{V}_{\lambda}^+) 
\cong 
\operatorname{Ind}_{\frak{A}_{n}(q)}^{\mathcal{H}_{n}(q)}(\Tilde{V}_{\lambda}^-) 
\cong V_{\lambda} 
\]
\end{Prop}
\begin{proof}
Directly from Proposition 6.6 we have
\[
\operatorname{Hom}_{\mathcal{H}_{n}(q)}
(\operatorname{Ind}_{\frak{A}_{n}(q)}^{\mathcal{H}_{n}(q)}(\Tilde{V}_{\lambda}),V_{\mu})
\cong 
\operatorname{Hom}_{\frak{A}_{n}(q)}
(\Tilde{V}_{\lambda},\operatorname{Res}_{\frak{A}_{n}(q)}^{\mathcal{H}_{n}(q)}(V_{\mu})) 
\]
and the Schur's Lemma means that its dimension is the multiplicity of $\Tilde{V}_{\lambda}$ 
in the restriction 
$\operatorname{Res}_{\mathcal{H}_{n}(q)}^{\frak{A}_{n}(q)}(V_{\mu})$. 
Multiplicities of irreducible $\frak{A}_{n}(q)$-modules in the restrictions of 
irreducible $\mathcal{H}_{n}(q)$-modules are already shown in Proposition 6.2 and Theorem 6.5 
for the non self-conjugate case and in Proposition 6.3 and Theorem 6.5 for the self-conjugate case.
\end{proof}

\end{document}